%&amstex
\magnification=\magstep1
\input amstex
\NoBlackBoxes
\documentstyle{amsppt}

\output={\plainoutput} 
\nopagenumbers 
\headline={\ifnum\pageno=1\hfil\else\ifodd\pageno\rhdg\else\lhdg\fi\fi} 
\def\rhdg{\tenrm\hfil GELFAND-DIKII FLOWS\hfil\folio} 
\def\lhdg{\tenrm\folio\hfil RICHARD BEALS \& D.H. SATTINGER\hfil} 
\voffset=2\baselineskip 
 
\font\auf = cmcsc10 at 12pt 

\hsize=5in                          
\vsize=7in     

\parindent=18pt 
\hfill 
\vskip 0.125in 
\centerline{\auf Action-Angle Variables for the Gel'fand-Dikii Flows} 
\vskip 0.375in 
\centerline{\auf Richard Beals \& D.H. Sattinger\footnote{Research supported
in part by 
N.S.F. Grants DMS-8916968 and DMS-8901607}}
\vskip 0.2in 
\centerline{Yale University and the University of Minnesota}
\vskip 0.2in 
\centerline{\it Dedicated to Professor Klaus Kirchg\"assner}
\vskip 0.375in
{\eightpoint \narrower 
Using the scattering transform for $n^{th}$ order linear scalar operators,
the Poisson bracket found by Gel'fand and Dikii, which generalizes the
Gardner Poisson bracket for the KdV hierarchy, is computed on the scattering
side.  Action-angle variables are then constructed.  Using this, complete
integrability is demonstrated in the strong sense.  Real action-angle
variables are constructed in the self-adjoint case.\par } 
\vskip .3cm 
\centerline{ Z. Angew. Math. Phys. {\bf 43} (1992) 219-242.}
\vskip .3cm
\noindent {\bf 1. The Gel'fand-Dikii hierarchies.}
\vskip .3cm
In this paper we construct action-angle variables and show complete
integrability in a strong sense for the Gel'fand-Dikii flows, which
generalize the KdV hierarchy.  Action-angle variables were first constructed
for the KdV hierarchy in 1971 by Zakharov and Faddeev [ZF]; and for the
nonlinear Schr\"odinger hierarchy in 1974 by Zakharov and Manakov [ZM]. 
McLaughlin [Mc] derived action-angle variables for the Toda lattice and
Sine-Gordon equations.

In a previous paper [BS1] we constructed action-angle variables and showed
complete integrability for a general class of Hamiltonian hierarchies based
on first order $n\times n$ isospectral operators.  It was shown that in some
cases, e.g. the three wave interaction, it was necessary to introduce flows
{\it nonlinear} in the scattering data in order to obtain the maximal number
of commuting Hamiltonian flows required for complete integrability.  This is
in contrast to the KdV and nonlinear Schr\"odinger equations in which 
the necessary commuting flows can be taken to be linear in the scattering data. 
 Consequently, the method of Lax  
pairs does not in general yield enough commuting flows for complete integrability.  In
the case of the Gel'fand-Dikii flows, however, the Lax pairs themselves do yield
enough commuting flows for complete integrability, and these flows are
linear in the scattering data.
 
Gel'fand and Dikii [GD] constructed a hierarchy of Hamiltonian flows
which generalize the KdV hierarchy based on the Schr\"odinger operator
$L=D^2+u$. They considered flows based on 
the $n^{th}$ order scalar differential operator
$$
L = \sum^n_{j=0} u_j(x)D^j;\quad D = {1\over i}{d\over dx};\quad u_n =
1,\quad u_{n-1} = 0 .
$$
where $u_j = u_j(x),~j < n-1 ,$ are elements of the Schwartz class of functions ${\Cal
S}(R)$.  The corresponding flows are given by
$$
{dL\over dt} = [L^{k/n}_+,L],\quad k=1,2,\dots
\eqno(1.1)
$$
where $k/n$ is not an integer, $L^{1/n}$ denotes the formal $n^{th}$ root of
$L$ considered as a pseudodifferential operator with symbol
$$
b(x,\xi) = \sum_{j\leq 1} w_j(x)\xi^j ;
$$
and $B^k_+$ denotes the differential part of a pseudodifferential operator
$B^k$. In the case $n=2$, $L$ is the Schr\"odinger operator, and (1.2) is the 
well-known Korteweg-deVries hierarchy (k odd).  For $n=3$ and $k=2$ equations
(1.1) yield the Boussinesq equation.
 
Gardner [G] showed that the KdV equation is Hamiltonian with respect to a
natural Poisson bracket.  Gel'fand and Dikii showed that the flows (1.1) are
Hamiltonian with respect to the Poisson bracket
$$
\{F,G\} = i\int^\infty_{-\infty} \sum_{r,s} [\ell_{rs} {\delta F\over\delta
u_s}] {\delta G\over\delta u_r}dx
\eqno(1.2)
$$
for the $n^{th}$ order case, where $\ell$ is the matrix of differential
operators
$$
\ell_{rs} = \sum^{n-1-r-s}_{k=0} [\bigl( {k+r\atop r}\bigr) u_{r+s+k+1}D^k -
\bigl( {k+s\atop s}\bigr) (-D)^ku_{r+s+k+1}]
\eqno(1.3)
$$
(We have introduced the factor $i$ in (1.2), by comparison with [GD], so that
the bracket is real in the self-adjoint case $L = L^\ast $.)
We assume throughout that  the $u_j, j<n-1 ,$ 
belong to the Schwartz class ${\Cal S}(R)$, so the terms $(-
D)^ku_{k+r+s+1}$ can be integrated by parts. After some
manipulation, (1.2) can be written
$$
\{F,G\} =
$$
$$
\int^\infty_{-\infty}\sum^n_{\gamma=1}\sum^{\gamma-1}_{r=0}\sum^r_{s=0}
u_\gamma\bigl( {r\atop s}\bigr) [(D^{r-s}{\delta F\over\delta u_{\gamma -r-
1}}){\delta G\over\delta u_s} - {\delta F\over\delta u_s} D^{r-s} {\delta
G\over\delta u_{\gamma -r-1}}]dx
\eqno(1.4)
$$
We will need this form below.
 
The scattering transformation of (1.1) as developed in [BDT] is based on the
wave functions $\psi$ which satisfy $L\psi = z^n\psi$.  This $n^{th}$ order
scalar equation can be written as a first order system in the standard way
 as follows.  Let
$e_{jk}$ denote the matrix with a 1 in the $jk$ place and zeroes elsewhere. 
Let $\psi_1$ be a solution of the scalar equation; the column vector $\psi =
(\psi_1,\psi_2\dots\psi_n)^t$ satisfies the first order system $(D-J_z-q)\psi
= 0$, where
$$
J_z =\left[ \matrix 0 & 1 & 0 &\hdots & 0\\
0 & 0 & 1 &\hdots & 0\\
& &\hdots & & &\\ 
0 & 0 & 0 &\hdots & 1\\
z^n & 0 & 0 &\hdots & 0\endmatrix \right]
\quad\quad q = - \left[ \matrix 0 & 0 &\hdots & 0 & 0\\
0 & 0 &\hdots & 0 & 0\\
& & &\hdots & &\\
0 & 0 &\hdots & 0 & 0\\
u_0 & u_1 & \hdots & u_{n-1} & 0\endmatrix \right]
$$
We denote by $P$ the linear space of all such matrices $q$ with entries in
${\Cal S}(R)$.  Let $\Phi$ and $\Psi$ be $n\times n$ matrices which satisfy
the first order systems
$$
D\Phi = (J_z +q)\Phi
\eqno(1.5a)
$$
$$
D\Psi = -\Psi (J_z+q)
\eqno(1.5b)
$$
These wave functions are uniquely determined by specifying their
normalization as $x\to -\infty$; that normalization is given in (1.9) below.
Define the matrix $a$ by $a=\Psi \Phi$; it is easily seen from (1.5) that $a$
is independent of $x$.

The entries of the first row of $\Phi$ satisfy the equations $L\phi =
z^n\phi$ while those of the last column of $\Psi$ satisfy the equation
$L^\dag\Psi = z^n\Psi$, where $L^\dag$ is the transpose of the operator $L$:
$$
L^\dag\Psi = \sum^n_{j=0} (-D)^ju_j\Psi
$$
Moreover, (1.5a) shows that $\Phi$ is the Wronskian of the functions of its
first row.

For $n\ge3$ let $\Sigma$ be the set of rays in the complex plane given by
$$
\Sigma = \{z\in C:Re\ i\alpha^jz = Re\ i\alpha^kz\ \text{for\ some}\ 0\leq
j<k<n\}
$$
where $\alpha = \exp 2\pi i/n$ is the primitive $n^{th}$ root of unity.  We
label the rays $\Sigma_j,\ j=0,1,\dots 2n-1$, with $\Sigma_0$ the negative
imaginary axis, and the numbering proceeding in the counterclockwise
direction.  We label domains of the complement
$C\backslash\Sigma$ by $\Omega_0,\Omega_1,\dots$ . (When $n=2$ the appropriate
set is $\Sigma = R$ and minor bookkeeping changes must be made in the 
following analysis; since the results were known for $n=2$ we omit the details.) 
\vfill\eject

\vglue3truein 

\centerline{Fig. 1. Fundamental Domain}
\vskip .2cm
\noindent It will sometimes be convenient, as in [BDT], to work in a {\it
local ordering}, that is, in an ordering that depends on $z$.  We define the
$j^{th}$ ordering of the roots to be $\{\alpha_1,\alpha_2\dots\alpha_n\}$,
where
$$
Re\ i\alpha_1z >Re\ i\alpha_2z >\dots > Re\ i\alpha_nz \quad\quad
\text{for }z\in\Omega_j .
\eqno(1.6)
$$
The ordering in $\Omega_0$ is given by $\{1,\alpha,\alpha^{-
1},\alpha^2\dots\}$ while that in $\Omega_1$ is given by  $\{1,\alpha^{-
1},\alpha,\alpha^{-2},\dots\}$.  The rays $\Sigma_0$ and $\Sigma_1$ are given
by
$$
\Sigma_0 = \{\xi :Im\ \xi<0,\ Re\ i\alpha^j\xi = Re\ i\alpha^{-j}\xi,\ j = 1,\dots,
[{n\over 2}]\},
$$
$$
\Sigma_1 = \{\xi :Im\ \xi<0,\ Re\ i\ \alpha^{-j}\xi = Re\ i\ \alpha^{j+1}\xi,\ j =
0,\dots,[{n\over 2}] -1\}.
$$
The transformation of a matrix from the $j^{th}$ to the $(j+1)^{st}$ ordering
is accomplished by conjugation by a permutation matrix $\pi_j$.  We see that
$\pi_j = \pi_{j+2}$, and that the structure of the permutations is different
for even and odd $n$.  The associated permutations are
$$
\pi_0 = (23)\dots (n-1,n)\quad \text{and}\quad \pi_1 = (12)\dots (n-2,n-1),\quad 
n\ \text{odd};
$$
and
$$
\pi_0 = (23)\dots (n-2,n-1)\quad \text{and}\quad \pi_1 = (12)\dots (n-1,n),\quad n\
\text{even}.
$$
 
We define the diagonal matrix $J(z)$ or $J_j$ to be the diagonal matrix with
entries $(\alpha_1,\alpha_2,\dots\alpha_n)$ where the roots are in the
$j^{th}$ ordering.  Similarly, let
$$
\Lambda_j = \left[ \matrix 1 & 1 & \hdots & 1\\
\alpha_1 & \alpha_2 & \hdots & \alpha_n\\
& & \hdots & \\  
\alpha^{n-1}_1 & \alpha^{n-1}_2 & \hdots & \alpha^{n-1}_n \endmatrix \right]
$$
where the roots are in the $j^{th}$ ordering; and for $z\in\Omega_j$, define
$\Lambda_z = d(z)\Lambda_j$, where $d(z) = diag (1,z,z^2,\dots z^{n-1})$. 
Then $\Lambda_j\Lambda_{j }^{\ast}= nI$ where $I$ is the $n\times n$ identity matrix,
and $J_z\Lambda_z = z\Lambda_zJ(z)$ for each $z\in C\backslash\Sigma$.
 
We now make the gauge transformations $\Phi = \Lambda_z\phi ,\Psi =
\psi\Lambda^{-1}_z$, $q\Lambda_z = \Lambda_zq_z$ to get the equations $(J =
J(z))$
$$
D\phi = (zJ+q_z)\phi,
\eqno(1.7a)
$$
$$
D\Psi = -\Psi (zJ+q_z).
\eqno(1.7b)
$$
Recall that $a(z) = \Psi\Phi$. It is clear that $a(z)$ is invariant under the
gauge transformation $\varphi\to\Phi,\ \psi\to\Psi$; hence $a(z) = \psi
(x,z)\phi (x,z)$.  In [BDT] fundamental solutions of systems (1.7a,b) are
constructed in each component of $C\backslash\Sigma$ by a technique involving
Volterra integral equations for successive exterior products of their columns.  The
solutions are constructed to satisfy the asymptotic conditions
$$
\lim_{x\to -\infty} \phi e^{iszJ} = 1,\quad\quad \limsup_{x\to +\infty} \Vert
\phi e^{-ixzJ}\Vert <+\infty;
\eqno(1.8a)
$$
$$
\lim_{x\to\infty} e^{ixzJ} \Psi = 1,\quad\quad \limsup_{x\to -\infty}\Vert
e^{ixzJ}\Psi\Vert <+\infty,
\eqno(1.8b)
$$
where $J = J(z)$.  The fundamental solutions $\Phi$ and $\Psi$ satisfy the
asymptotic conditions
$$
\lim_{x\to -\infty} \Lambda^{-1}_z \Phi e^{-ixzJ(z)} = 1,\quad\quad
\lim_{x\to\infty}e^{ixzJ(z)}\Psi\Lambda_z = 1.
\eqno(1.9)
$$
As $x\to -\infty$, respectively $+\infty,\Phi$ and $\Psi$ are asymptotic to
$U =\Lambda_z e^{ixzJ(z)}$ and $U^{-1}$ respectively; a
simple calculation shows that
$$
U_{jk} = D^{j-1}e^{ix\alpha_kz},\quad\quad (U^{-1})_{jk} = (-D)^{n-k}
{\alpha_j\over nz^{n-1}} e^{-ixz\alpha_j}.
\eqno(1.10)
$$
 
Let $\Sigma_j$ be a ray of $\Sigma$ and let $\xi$ belong to $\Sigma_j$.  Denote by
$\phi_\pm$ the boundary values of the wave function $\phi$ defined on the
components of $C\backslash \Sigma$ separated by $\Sigma_j$.  In [BDT] the
{\it scattering data} for $L$ associated with $\Sigma_j$ is defined by the jump
relations $v(\xi) = \phi^{-1}_-\phi_+\pi_j$.  In general, the scattering data consists of
discrete data, as well, which lives on those points $z$ for which the
construction of $\phi$ and $\psi$ breaks down. We assume in this article
that there is no discrete component of the scattering data, so that $L$ is
uniquely determined by its continuous scattering data.
 
As in [BC] and [BS1] we may instead define a {\it scattering matrix} for the first order
system (1.7a) as follows.  Fundamental solution matrices which satisfy the asymptotic
conditions (1.8) and (1.9) can also be
constructed for $\xi\in\Sigma\backslash 0$. 
For
$\xi\in\Sigma$, let $\Pi_\xi$ denote the projection
$$
(\Pi_\xi s)_{jk} = \cases s_{jk} &\text{if}\ Re\ i\xi (\lambda_j-\lambda_k) =
0\cr
0\ &\text{otherwise}.\endcases
$$
Then the limit
$$
s(\xi) = \lim_{x\to +\infty} \Pi_\xi e^{-ix\xi J}\phi (x,\xi) =
\lim_{x\to+\infty} \Pi_\xi e^{-ix\xi J}\Lambda^{-1}_\xi \Phi (x,\xi)
\eqno(1.12)
$$
exists.  The matrices $s(\xi)$, as $\xi$ ranges over the rays $\Sigma_j$, 
determine the continuous data  $v$, and, in the absence of discrete data,
$v$ determines $s$. In \S 3
action-angle variables for the Gel'fand-Dikii flows are constructed directly
from the entries of $s$, not $v$, so it is important to pass to  $s$.
\vskip .3cm
\noindent {\bf Lemma 1.1.} $s = a = \psi\phi$.
\vskip .3cm
\noindent Proof: Writing $\phi = me^{ix\xi J}$ and $\psi = e^{-ix\xi J}w$,
the functions $m$ and $w$ are bounded on $-\infty <x<\infty$ by (1.8) and
(1.9).  Therefore $wm = e^{ix\xi J}a\ e^{-ix\xi J}$ is bounded for $x\in R$,
and so $\Pi_\xi a = a$.  Now
$$
\eqalign{s(\xi) &= \lim_{x\to\infty}\Pi_\xi e^{-ix\xi J}\phi =
\lim_{x\to\infty} \Pi_\xi e^{-ix\xi J}\psi^{-1}a\cr
&= \lim_{x\to\infty}\Pi_\xi e^{-ix\xi J}W^{-1}e^{ix\xi J}a =
\lim_{x\to\infty} [\Pi_\xi e^{-ix\xi J}w^{-1}e^{ix\xi J}]a\cr}
$$
since $\Pi_\xi (Ba) = (\Pi_\xi B)a$ for any matrix $B$ whenever $a$ is in the
range of $\Pi_\xi$.  As $x$ tends to $\infty$, $w^{-1}$ tends I; and since
$\Pi_\xi$ projects onto the oscillatory part of $\exp\{ix\xi adJ(\xi)\}w$,
$$
\lim_{x\to\infty} \Pi_\xi e^{-ix\xi J} w^{-1}e^{ix\xi J} = I .
$$
Hence $s=a$ and Lemma 1.1 is proved.
\vskip .3cm 
Since $a$ lies in the range of the projection $\Pi_\xi$ its only nonzero
entries occur where $Re\ i\alpha_j\xi = Re\ i\alpha_k\xi$.  In the local
ordering, $a(\xi)$ has a block diagonal structure consisting of diagonal
entries or $2\times 2$ diagonal blocks. We label the $2\times 2$ blocks by
 $B_\nu,\nu = 1,\dots n-1$ as indicated. The block diagonal structure of
$a(\xi)$ is different for even and odd $n$ as indicated in the following
diagrams:
                     
\newpage 

\centerline{{\bf n odd}}
\vskip-.7truein
$$
\pmatrix
a_{11}
\\
&\boxed{B_2}\\
\\
&&\boxed{B_4}\\
&&&\ddots\\
\\
&&&&\boxed{\!B_{n-1}\!}
\endpmatrix
\matrix
\\
\\
\\
\\
\\
\\                          
\\
\\
\\
\\
\\
\\
\!\!\!\!\!\!\!\!\!\!\!\!\!\!\!\!\!\!\!\!\!
\!\!\!\!\!\!\!\!\!\!\!\!\!\!\!\!\!\!\!\!\!\!\!\!\!\!\!\!\!\!\!\!\!\!\!\!\!\!\!\!\!\!\!
\!\!\!\!\!\!\!\!\!\!\!\!\!\!\!\!\!\!\!\!\!\!\!\!\!\!\!\!\!\!\!\!\!\!\!\!\!\!\!\!\!\!\!\!\!\!
\Sigma_0&&&&&&&&&&&&
\endmatrix
\!\!\!\!\!\!\!\!\!\!\!\!\!\!\!\!\!\!\!\!\!\!\!\!\!\!\!\!\!\!\!\!\!\!\!\!\!\!\!\!\!\!\!
\!\!\!\!\!
\pmatrix
\boxed{B_1}\\ 
\\
&\boxed{B_3}\\
\\
&&\ddots\\
\\
&&&a_{nn}
\endpmatrix          
\matrix
\\
\\
\\
\\
\\
\\                          
\\
\\
\\
\\
\\
\\
\!\!\!\!\!\!\!\!\!\!\!\!\!\!\!\!\!\!\!\!\!\!\!\!\!\!\!\!\!\!\!\!\!\!\!\!\!
\!\!\!\!\!\!\!\!\!\!\!\!\!\!\!\!\!\!\!\!\!\!\!\!\!\!\!\!\!\!\!\!\!\!\!\!\!\!\!\!\!\!\!\!\!\!
\Sigma_1&&&&&&&&&&&&
\endmatrix
$$        
\medskip
\centerline{{\bf n even}}
\vskip-.7truein
$$
\pmatrix
a_{11}\\
\\
&\boxed{B_2}\\
\\         
&&\ddots\\
&&&\boxed{B_{n-2}}\\
\\               
&&&&a_{nn}
\endpmatrix
\matrix
\\
\\
\\
\\
\\
\\                          
\\
\\
\\
\\
\\
\\
\!\!\!\!\!\!\!\!\!\!\!\!\!\!\!\!\!\!\!\!\!
\!\!\!\!\!\!\!\!\!\!\!\!\!\!\!\!\!\!\!\!\!\!\!\!\!\!\!\!\!\!\!\!\!\!\!\!\!\!\!\!\!\!\!\!\!\!
\!\!\!\!\!\!\!\!\!\!\!\!\!\!\!\!\!\!\!\!\!\!\!\!\!\!\!\!\!\!\!\!\!\!\!\!\!\!\!\!\!\!\!\!\!\!
\Sigma_0&&&&&&&&&&&&
\endmatrix
\!\!\!\!\!\!\!\!\!\!\!\!\!\!\!\!\!\!\!\!\!\!\!\!\!\!\!\!\!\!\!\!\!\!\!\!\!\!\!\!\!\!\!
\!\!\!\!\!
\pmatrix
\boxed{B_1}\\
\\
&\ddots\\
&&&\ddots\\
&&&&\boxed{B_{n-1}}\\
\endpmatrix
\matrix
\\
\\
\\
\\
\\
\\                          
\\
\\
\\
\\
\\
\\
\!\!\!\!\!\!\!\!\!\!\!\!\!\!\!\!\!\!\!\!\!\!\!\!\!\!\!\!\!\!\!
\!\!\!\!\!\!\!\!\!\!\!\!\!\!\!\!\!\!\!\!\!\!\!\!\!\!\!\!\!\!\!\!\!\!\!\!\!\!\!\!\!\!\!\!\!\!
\Sigma_1&&&&&&&&&&&&
\endmatrix
$$     
$$
B_\nu=\pmatrix
a_{\nu \nu }&&a_{\nu,\nu+1}\\ \\
a_{\nu+1, \nu} &&a_{\nu+1, \nu+1}
\endpmatrix\qquad \qquad \Delta_\nu =\det B_\nu
$$                  
\medskip
\centerline{ Fig. 2 Block Diagonal Structure}
\medskip
\noindent{\bf Lemma 1.2}.\quad {\it The wave functions $\Phi$ and $\Psi$ and
scattering data $a(\xi)$ are invariant under the rotations $z\to\alpha z$; in
particular, $a(\xi) = a(\alpha\xi)$.  Consequently the scattering data is
determined entirely by the matrices $a(\xi)$ for $\xi$ on the rays $\Sigma_0$
and $\Sigma_1$.}
\vskip .3cm
In Lemma 1.2 it is understood that we are working in the local
representation.  The invariance of the wave functions and scattering data
under rotation by $\alpha$ means that all data is contained in the wave
functions defined on the fundamental domain in Fig. 1. 
\vskip .3cm
\noindent Proof:\quad First note that $J_z,\Lambda_z$, and $zJ(z)$ are
invariant under the rotation $z\to\alpha z$, since the matrix representations
are local.  Consequently, equations (1.5a,b) and the matrix $U$ of (1.10) are
invariant.  By uniqueness, the wave functions $\Phi$ and $\Psi$ are
invariant, and therefore so too is $a(\xi) = \Psi (x,\xi)\Phi (x,\xi)$.
\vskip .3cm
 
Later we show that the scattering data possesses an additional symmetry when
$L$ is self-adjoint.      

\vskip .3cm
\noindent{\bf 2. The Poisson brackets.}
\vskip .3cm
We now turn to the computation of the Poisson brackets of 
the scattering variables.  The entries of the scattering matrix $a(\xi)$ are
regarded as functionals on $P$, and we compute their corresponding Poisson
brackets.  As in [BS1], the gradients $\delta a_{jk}(\xi)/\delta q$ do not
decay as $\vert x\vert\to\infty$. Therefore the integral in (1.2) does not have an
absolutely convergent integrand, and it is necessary instead to use the
regularization
$$
\lim_{N\to\infty} i\int^N_{-N}\sum_{r,s} [\ell_{rs} {\delta
a_{jk}(\xi)\over\delta u_s}] {\delta a_{lm}(\eta)\over\delta u_r} dx.
$$
This limit exists only in the sense of distributions; so the precise
interpretation of the calculation is as follows.  Given any pair of functions
$f,g\in {\Cal D}(\Sigma)$, define the functionals 
$$
F(q) = \int_{\Sigma} a_{jk}(\xi)f(\xi)d\xi,\quad\quad G(q) = \int_{\Sigma}
a_{l m}(\eta)g(\eta)d\eta.
$$
Then the formal calculation
$$
\eqalign{\{F,G\} &= \lim_{N\to\infty} i\int^N_{-N} \sum_{r,s}[ \ell_{rs}
{\delta F\over\delta u_s}] {\delta G\over\delta u_r} dx\cr
&= \lim_{N\to\infty} i\int_{\Sigma}\int_{\Sigma}\int^N_{-N} \sum_{r,s}
[\ell_{r,s} {\delta a_{jk}(\xi)\over\delta u_s}] {\delta
a_{l m}(\eta)\over\delta u_r} f(\xi) g(\eta)dx d\xi d\eta\cr
&= \int_{\Sigma}\int_{\Sigma} \{a_{jk}(\xi),a_{l m}(\eta)\}
f(\xi)g(\eta)d\xi d\eta\cr}
$$
furnishes the defining relation for the distribution
$\{a_{jk}(\xi),a_{l m}(\eta)\}\in {\Cal
D}'(\Sigma\times\Sigma)$.
 
Let $\delta\Phi$ and $\delta a$ denote variations in the wave function and
scattering data due to a variation in the potential matrix $q$ (i.e. in the
coefficients of $L$).  The following identity may be established
by replacing $q$ by $q+\epsilon\delta q$, forming the equation for the
difference in the wave functions, and taking the limit as $\epsilon\to 0$
(cf. \S D of [BC]):
$$
D(\Psi\delta\Phi) = \Psi\delta q\Phi .
\eqno(2.1)
$$
Integrating this identity over $(-\infty,\infty)$ and using the fact that
$\delta\phi$ and $\delta\psi$ vanish at $-\infty$ and $+\infty$ respectively,
we obtain
$$
\delta a =i \int^\infty_{-\infty} \Psi\delta q\Phi dx .
$$
Therefore
$$
(-i)\delta a_{jk} = \int^\infty_{-\infty} \text{tr}[\Psi\delta q\Phi e_{kj}]dx =
\sum^{n-1}_{l =0} \int^\infty_{-\infty} \text{tr}[\Phi e_{kj}\Psi e_{n,l
+1}]\delta u_l dx.
$$
Hence
$$
(-i){\delta a_{jk}(\xi)\over\delta u_l} = (\Phi e_{kj}\Psi)_{l +1,n} =
\Phi_{l +1,k}(x,\xi)\Psi_{j,n}(x,\xi) = (D^l\phi_k)\psi_j
\eqno(2.2)
$$
where $\phi_k$ and $\psi_j$ are the entries of the first row of $\Phi$ and
last column of $\Psi$.  The Poisson bracket between entries of the scattering
matrix $a$ is therefore (using (1.4))
$$
\eqalign{-i\int^\infty_{-\infty} \sum^n_{\gamma=1}\sum^{\gamma-
1}_{r=0}\sum^r_{s=0} u_\gamma \bigl( {r\atop s}\bigr) &\{[D^{r-
s}(\psi_jD^{\gamma -r-1}\phi_k)[\psi_l D^s\phi_m]\cr 
&- [\psi_jD^s\phi_k][D^{r-s}\psi_l D^{\gamma -r-1}\phi_m]\}dx.\cr}
\eqno(2.3)
$$
Here $\psi_j$ and $\phi_k$ are evaluated at $(x,\xi)$ and $\psi_l$ and
$\phi_m$ are evaluated at $(x,\eta)$.
 
When $\xi$ and $\eta$ are on different rays we must choose a common
representation for the matrices $a(\xi)$ and $a(\eta)$; in the following, we
work in a fixed representation.
 
The integral in (2.3) is evaluated by expressing the integrand as an exact
derivative.  More precisely,
\vskip .4cm
\noindent{\bf Lemma 2.1}\quad {\it Let $u_\gamma, \gamma < n-1,$ belong to the Schwartz
class.  Then
$$
\{a_{jk}(\xi),a_{lm}(\eta)\} = {-i\over\xi^n-\eta^n} \int^\infty_{-\infty}
D[W(\phi_k(\xi),\psi_l(\eta))W(\phi_m(\eta),\psi_j(\xi))]dx,
\eqno(2.4)
$$
where}
$$
W(f,g) = \sum^n_{\gamma=1} \sum^{\gamma-1}_{r=0} (D^{\gamma-r-1}f)(-D)^r(u_\gamma
g).
\eqno(2.5)
$$
 
\noindent Proof:  For scalar functions $f$ and $g$,
$$
gLf-fL^\dag g = DW(f,g) .
$$
In particular, if $Lf = \xi^nf$ and $L^\dag g = \eta^ng$, then
$$
DW(f,g) = (\xi^n-\eta^n)fg.
\eqno(2.6)
$$
It follows from (2.5) and (2.6) that
$$
{1\over\xi^n-\eta^n}
D[W(\phi_k(\xi),\psi_l(\eta))W(\phi_m(\eta),\psi_l(\xi))]
$$
$$
= \phi_k\psi_l W(\phi_m,\psi_j)-\phi_m\psi_j W(\phi_k,\psi_l)
$$
$$
= \sum^n_{\gamma=1} \sum^{\gamma-1}_{r=0} \{\phi_k\psi_l [(-D)^r
(u_\gamma\psi_j)]D^{\gamma-r-1}\phi_m-\phi_m\psi_j[(-D)^r u_\gamma \psi_l
]D^{\gamma-r-1}\phi_k\}.
$$
Since we assume the $u_\gamma,~\gamma < n-1~,$ are in ${\Cal S}(R)$ this 
expression is equivalent modulo an exact
derivative to
$$
\sum^n_{\gamma=1} u_\gamma \sum^{\gamma-1}_{r=0} \psi_jD^r (\phi_k\psi_l
D^{\gamma-r-1}\phi_m)-\psi_l D^r (\phi_m\psi_jD^{\gamma-r-1}\phi_k).
$$
Comparing this with (2.3) we get (2.4).
 
We now evaluate the integral in (2.4).  As noted in the beginning of this
section the integral in (2.4) must be interpreted in the sense of
distributions.  We obtain
$$
\{a_{jk}(\xi),a_{lm}(\eta)\} =- \lim_{N\to\infty} {W(\phi_k(\xi),\psi_l
(\eta))W(\phi_m(\eta),\psi_j(\xi))\over\xi^n-\eta^n} \Bigr|^N_{-N}
$$
where $\phi_k(\xi)$ denotes $\phi_k(x,\xi)$, etc.  Since $u_n = 1$ and the
remaining coefficients $u_j$ tend to zero rapidly as $x\to\pm\infty$,
$$
W(\phi,\psi)\sim\ W_0(\phi,\psi) = \sum^{n-1}_{s=0} (D^{n-s-1}\phi)(-D)^s\psi
\quad \text{as}\ \vert x\vert\to\infty.
$$
 
From $\Psi\Phi = a$ and (1.10) we have, as  $x\to -\infty$ and $x\to +\infty$
respectively,
$$
\phi_j(x,\xi)\sim\ e^{ix\alpha_j\xi}\quad\quad \phi_j(x,\xi)\sim \ \sum_r
e^{ix\xi\alpha_r} a_{rj}(\xi)
$$
and
$$
\psi_j\sim\sum_r a_{jr}(\xi) {\alpha_r\over n\xi^{n-1}} e^{-
ix\xi\alpha_r}\quad\quad \psi_j\sim {\alpha_j\over n\xi^{n-1}} e^{-
ix\xi\alpha_j} .
$$
From the block diagonal nature of $a(\xi)$, the sums for $\phi_j$ and
$\psi_j$ extend only over the one or two element set $\{r:\text{Im}\ \xi(\alpha_r-
\alpha_j)=0\}$.
 
Straightforward calculations show that
$$
W_0(e^{ix\xi\alpha_k},e^{-ix\eta\alpha_p}) = i^{n-1}{\xi^n-
\eta^n\over\alpha_k\xi-\alpha_p\eta} e^{ix(\alpha_k\xi -\alpha_p\eta)}
\eqno(2.7)
$$
and that as $x\to\pm\infty$
$$
{1\over\xi^n-\eta^n}
W(\phi_k(\xi),\psi_l(\eta))W(\phi_m(\eta),\psi_j(\xi))
$$
is asymptotic to
$$
{(-1)^n(\xi^n-\eta^n)\over n^2(\xi\eta)^{n-1}} \sum_{p,r}\alpha_j\alpha_l
a_{pk}(\xi)a_{rm}(\eta) {e^{ix(\xi\alpha_p-\eta\alpha_l + \eta\alpha_r
-\xi\alpha_j)}\over (\xi\alpha_p-\eta\alpha_l)(\eta\alpha_r-\xi\alpha_j)}
\eqno(2.8+)
$$
and
$$
{(-1)^n(\xi^n-\eta^n)\over n^2(\xi\eta)^{n-1}} \sum_{p,r}\alpha_p\alpha_r
a_{lp}(\eta)a_{jr}(\xi) {e^{ix(\xi\alpha_k-\eta\alpha_p+\eta\alpha_m-
\xi\alpha_r)}\over (\xi\alpha_k-\eta\alpha_p)(\eta\alpha_m-\xi\alpha_r)} 
\eqno(2.8-)
$$
respectively. Hence
$$
\{a_{jk}(\xi),a_{lm}(\eta)\} = - \sum_{p,r} a_{pk}(\xi)
a_{rm}(\eta)R^+_{pl,rj}(\xi,\eta)-a_{lp}(\eta)a_{jr}(\xi)R^-
_{kp,mr}(\xi,\eta)
\eqno(2.9)
$$
where
$$
R^{\pm}_{pl,rj}(\xi,\eta) = {(-1)^n(\xi^n-\eta^n)\over n^2(\xi\eta)^{n-1}}
\alpha_l\alpha_j \lim_{x\to\pm\infty} {e^{ix(\xi(\alpha_p-
\alpha_j)+\eta(\alpha_r-\alpha_l)}\over (\xi\alpha_p-
\eta\alpha_l)(\eta\alpha_r-\xi\alpha_j)} \ .
\eqno(2.10)
$$
 
\noindent {\bf Lemma 2.2.} {\it For} $\xi,\ \eta$ {\it in} $\Sigma_0 \cup \Sigma_1$
{\it the limits in (2.10) are}
$$
\eqalign{R^\pm_{pl,rj}(\xi,\eta) &= {(-1)^{n-1}\over n^2(\xi\eta)^{n-
1}}\delta_{jp}\delta_{rl}\alpha_j\alpha_l {\xi^n-\eta^n\over\xi\alpha_j-
\eta\alpha_l} p.v. {1\over\xi\alpha_j-\eta\alpha_l}\cr
&\pm \pi i\delta_{lp}\delta_{jr} {(-1)^{n-1}\over n\xi^{n-1}}
\delta(\xi-\eta)sgn(l-j).\cr}
\eqno(2.11)
$$
Here $\delta (\xi -\eta )$ denotes the complex valued distribution defined with
respect to the complex measure $d\xi$ on the rays of $\Sigma$ such that
$$
\int_{\Sigma}f(\eta )\delta(\xi -\eta )d\eta =f(\xi ).
$$
Brackets for general $\xi, \eta$ are obtained
using the symmetry in Lemma 1.2.
\vskip .3cm
 
\noindent Proof:  We use the $0^{th}$ representation on
both rays.  Note that the block diagonal structure of  a  on  $\Sigma_1$ is
the same in either the $0^{th}$ or the $1^{st}$ representation and so is the
same as that shown in Fig. 1.2.
 
The first term in (2.11) arises when $\xi(\alpha_p-\alpha_j)+\eta (\alpha_r-
\alpha_l)$ vanishes identically in $\xi$ and $\eta$, that is when $j=p$ and
$r=l$.  In that case the argument of the exponential vanishes and we get a
simple principal value term.  Note that
$$
{\xi^n-\eta^n\over\xi\alpha_j-\eta\alpha_l}
\eqno(2.12)
$$
is regular, since $\alpha_j$ and $\alpha_l$ are roots of unity.
 
Due to the block structure of $a(\xi)$, the sums in (2.8) extend only over
the sets where $\text{Im}\ \xi\ (\alpha_p-\alpha_j) = \text{Im}\ \eta\ (\alpha_r-
\alpha_l) = 0$; hence the exponential term in (2.10) is purely oscillatory. 
We use the identity
$$
\lim_{x\to\pm\infty} p.v. {e^{i\lambda x(\xi-\eta)}\over\xi-\eta} = \pm\pi 
 i\ sgn(\lambda)\delta (\xi-\eta),\quad\quad \lambda\in R\backslash 0,\ (\xi-
\eta)\in R ,
\eqno(2.13)
$$
to evaluate the limit when the argument of the exponential does not vanish. 
Since (2.12) is regular, the limit in (2.10) vanishes by the Riemann-Lebesgue
lemma whenever only one of the factors $(\xi\alpha_p-
\eta\alpha_l),(\xi\alpha_j-\eta\alpha_r)$ vanishes.  Since $\Sigma_0$ and
$\Sigma_1$ are on different orbits under multiplication by roots of unity,
these factors can vanish only if $\xi$ and $\eta$ are on the same ray.  If
both factors vanish, then $j=r$ and $p=l$ and we must evaluate
$$
\lim_{x\to\pm\infty} {e^{ix(\xi-\eta)(\alpha_l-\alpha_j)}\over\xi-
\eta}
$$
when $Im(\xi-\eta)(\alpha_l-\alpha_j) = 0$.  In evaluating this limit, we
must keep in mind that it is defined in the sense of distributions; so we
really need to consider the limits
$$
\eqalign{&\lim_{x\to\pm\infty} \int_{\Sigma_r} {e^{ix(\xi-\eta)(\alpha_j-
\alpha_l)}\over (\xi-\eta)} f(\eta)d\eta \cr 
&= \lim_{x\to\pm\infty} \int_{\Sigma_r} {e^{ix(\xi-\eta)(\alpha_j-
\alpha_l)}\over (\xi-\eta)(\alpha_j-\alpha_l)} f(\eta)d(\alpha_j-
\alpha_l)\eta\cr}
$$
where $f$ is a smooth test function on the ray $\Sigma_r,\ r=0,1$.  Making
the change of variables $s=\xi(\alpha_j-\alpha_l),\ t = \eta(\alpha_j-
\alpha_l)$ we get
$$
\eqalign{&= \lim_{x\to\pm\infty} \int^{sgn(l-j)\infty}_0 {e^{ix(s-
t)}\over s-t} f({t\over\alpha_j-\alpha_l})dt\cr 
&= \pm \pi if({s\over\alpha_j-\alpha_l})sgn(l-j)\cr 
&= \pm\pi if(\xi)sgn(l-j).\cr}
$$
Therefore in the sense of distributions the above limit is
$$
\lim_{x\to\pm\infty} {e^{ix(\xi-\eta)(\alpha_l-\alpha_j)}\over\xi-\eta} =
\pm\pi i\delta(\xi-\eta)sgn(l-j).
$$      

Equation (2.11) follows, since
$$
{\xi^n-\eta^n\over\xi-\eta}\delta (\xi-\eta ) = n\xi^{n-1}\delta
(\xi-\eta ).
$$
 
Using lemma 2.2 and (2.10) we get
$$
\eqalign{\{&a_{jk}(\xi),a_{lm}(\eta)\} = a_{jk}(\xi)a_{lm}(\eta) {(-1)^{n}
\over n^2(\xi\eta)^{n-1}} [Q_{j,l}(\xi,\eta)-Q_{k,m}(\xi,\eta)]\cr
&+ (-1)^{n}\pi ia_{l k}(\xi)a_{jm}(\xi) {1\over n\xi^{n-1}}\delta
(\xi-\eta)[sgn(l -j)+sgn(k-m)]\cr}
\eqno(2.14)
$$
where
$$
Q_{j,l}(\xi,\eta) = {\xi^n-\eta^n\over\xi\alpha_j-\eta\alpha_l}
\alpha_j\alpha_l \  p.v. {1\over\xi\alpha_j-\eta\alpha_l}.
$$
 
In (2.14) all matrices are in a fixed (global) ordering.  We shall always
take the ordering to be that of $\Omega_0$.  For $n=2$ (2.14) gives the known
result for the KdV case.  (cf. [NMPZ] p. 41; the factor of 1/2 is due to
different normalizations of the Gardner Poisson brackets.)   

\vskip .3cm
\noindent {\bf 3. Action-angle variables; general case.}
\vskip .3cm         

The construction of action-angle variables for general first order $n\times
n$ isospectral operators was given in [BS1].  While the situation for the
$n^{th}$ order scalar case is similar, there are some differences.  The local
part of the bracket (2.14), i.e. the part containing the delta term, is the same,
apart from the factor $(n\xi^{n-1})^{-1}$, as in the case of
first order systems.  The nonlocal term, -- involving $p.v.(\xi\alpha_j-
\eta\alpha_l)^{-1}$ -- is more complicated than that in the
systems case.  Finally, because of the block structure of the scattering
data, we need only consider the case of $2\times 2$ matrices, for which the
structure of the action-angle variables is especially simple.   

It is convenient to use the positive real, rotationally invariant
 measure $\xi^{-1} d\xi$ on $\Sigma$. We denote by $\tilde \delta(\xi / \eta)$ 
the associated distribution defined by
$$
\int _{\Sigma} f(\xi) \tilde \delta (\xi / \eta ){d\xi \over \xi}=f(\eta).
$$      
A simple computation shows that $\tilde \delta(\xi / \eta )= \xi~ \delta 
(\xi - \eta)$, the latter being the distribution defined in the previous
section.
 
We label the $2\times 2$ blocks of the scattering matrix $a$ by $\nu =1,\dots
n-1$, as in Fig. 2, and define as the canonical variables associated with
the $\nu^{th}$ block,
$$
p_\nu (\xi) = n(-\xi)^n
\log {a_{\nu\nu}(\xi)a_{\nu+1,\nu+1}(\xi)\over\Delta_\nu},\quad\quad q_\nu(\xi) =
{1\over 4\pi i} \log {\alpha^\nu a_{\nu+1,\nu}(\xi)\over a_{\nu,\nu+1}(\xi)},
\eqno(3.1)
$$
where $\Delta_\nu(\xi) = a_{\nu\nu}a_{\nu+1,\nu+1}-
a_{\nu,\nu+1}a_{\nu+1,\nu}$ is the determinant of the $2\times 2$ block
$B_\nu$.  In (3.1), $\xi$ is on $\Sigma_0$ or $\Sigma_1$ according as $\nu$
is even or odd; and we work in the local representation.  The phase factor
$\alpha^\nu$ inserted in the definition of the $q_\nu$ is needed to get the
right symmetries in the self-adjoint case, discussed in the next section.  It
simply involves changing $q_\nu$ by an additive constant, and so does not
affect the Poisson brackets.

\vskip .3cm
\noindent{\bf Theorem 3.1.}\quad {\it The variables $p_\nu,q_\nu$ satisfy the
canonical commutation relations}
$$
\{p_\nu(\xi),p_\mu(\eta)\} = \{q_\nu(\xi),q_\mu(\eta\} = 0,\quad\quad
\{p_\nu(\xi),q_\mu(\eta)\} = \delta_{\mu \nu}~\tilde \delta (\xi / \eta )
$$
{\it with respect to the positive measure} $\xi^{-1}d\xi$. 
\vskip .3cm
\noindent Proof:\quad The distribution-valued Poisson bracket (2.15) is a sum
of two terms, a local term, involving the distribution $\delta (\xi -
\eta )$, and a nonlocal term, involving the principal value terms such
as p.v. $(\xi\alpha_j-\eta\alpha_l)^{-1}$.  (Note however, that
$(\xi\alpha_j-\eta\alpha_l)^{-1}$ is singular only if $\xi$ and $\eta$ are on
the same ray, since $\Sigma_0$ and $\Sigma_1$ are on different orbits under
multiplication by a root of unity.)  As in [BS1] we denote the local part of
the Poisson bracket by $[a_{jk},a_{lm}]$ and the nonlocal term by
$<a_{jk},a_{lm}>$.
 
Recall that the relations (2.14) refer all quantities to the $0^{th}$
representation, whereas the definition of the canonical variables in (3.1)
uses local representations on each of the rays.  To convert the
representation from $0^{th}$ to $1^{st}$, we conjugate by the permutation
$\pi_0$; and this interchanges $a_{\nu\nu}$ with $a_{\nu+1,\nu+1}$ and
$a_{\nu,\nu+1}$ with $a_{\nu+1,\nu}$; so $p_\nu$ is unchanged and $q_\nu$ is
replaced by $(const.-q_\nu)$.  Since the local component of the Poisson
bracket vanishes when $\xi$ and $\eta$ are on different rays, none of the
following computations are materially affected.
 
By (2.14) we have for all $j,k,\xi,\eta$ (with $j<k$, say)
$$
\{a_{jj}(\xi),a_{kk}(\eta)\} = a_{jj}(\xi)a_{kk}(\eta) {(-1)^n\over
n^2(\xi\eta)^{n-1}} [Q_{j,k}(\xi,\eta) - Q_{j,k}(\xi,\eta)]
$$
$$                                                              
+ (-1)^n\pi ia_{jk}(\xi)a_{kj}(\xi) {1\over n\xi^{n-1}} \delta (\xi-
\eta )[sgn(k-j)+sgn(j-k)] = 0 .
$$
The first set of relations in (3.1) follows if we can show that
$\{\Delta_\mu,\Delta_\nu\} = \{\Delta_\mu,a_{jj}\} = 0$ for any $\mu,\nu$,
and $j$.  Note that for the local part of the bracket,
$[a_{jk}(\xi),a_{lm}(\eta)] = 0$ whenever $\xi$ and $\eta$ are on different
rays, or whenever $a_{jk}$ is strictly ``northwest'' or ``southeast'' of
$a_{lm}$, due to the factor $[sgn(j-l)+sgn(m-k)]$. 
The same result can also be deduced from the block diagonal structure of $a$. 
Thus $[\Delta_\mu,\Delta_\nu] = [\Delta_\mu,a_{jj}] = 0$.  Let us consider
the nonlocal bracket $<\Delta_\nu,a_{jj}>$.  We have, by the derivation
property of the Poisson bracket (cf [BS1]),
$$
<\Delta_\nu,a_{jj}> = <a_{\nu\nu}a_{\nu+1,\nu+1} -
a_{\nu,\nu+1}a_{\nu+1,\nu},a_{jj}>
$$
$$
\eqalign{= a_{\nu\nu}<a_{\nu+1,\nu+1},a_{jj}> &+
a_{\nu+1,\nu+1}<a_{\nu\nu},a_{jj}>
-a_{\nu,\nu+1}<a_{\nu+1,\nu},a_{jj}>\cr
&- a_{\nu+1,\nu}<a_{\nu,\nu+1},a_{jj}>\cr}
$$
$$
= - {(-1)^n\over n^2(\xi\eta)^{n-1}}
a_{\nu,\nu+1}a_{\nu+1,\nu}a_{jj}[Q_{\nu+1,j}-Q_{\nu,j}+Q_{\nu,j}-Q_{\nu+1,j}]
= 0.
$$
A similar calculation shows that $<\Delta_\mu,\Delta_\nu> = 0$; and the first
set of commutation relations in (3.1) is proved.
 
Next,
$$
\{\log {a_{\mu+1,\mu}\over a_{\mu,\mu+1}},  \log {a_{\nu+1,\nu}\over
a_{\nu,\nu+1}}\}
= \{\log a_{\mu+1,\mu}-\log a_{\mu,\mu+1},\log a_{\nu+1,\nu}-\log
a_{\nu,\nu+1}\}
$$
$$
= {1\over a_{\mu+1,\mu}a_{\nu+1,\nu}} \{a_{\mu,\mu+1},a_{\nu+1,\nu}\} -
{1\over a_{\mu+1,\mu}a_{\nu,\nu+1}}\{a_{\nu+1,\nu},a_{\mu,\mu+1}\}
$$
$$
+{1\over a_{\mu,\mu+1}a_{\nu+1,\nu}}\{a_{\mu,\mu+1},a_{\nu+1,\nu}\} + {1\over
a_{\mu,\mu+1}a_{\nu,\nu+1}}\{a_{\mu,\mu+1},a_{\nu,\nu+1}\} .
$$
 
We claim that the local part of the above Poisson bracket vanishes.  It
certainly vanishes if $\xi$ and $\eta$ are on different rays.  If $\xi$ and
$\eta$ are on the same ray and $\mu <\nu$, the local component vanishes due
to the block structure of $a$; while if $\mu =\nu$ then the first and fourth
terms vanish while the second and third cancel.
 
Therefore, the above Poisson bracket reduces to the nonlocal terms, viz
$$
\eqalign{ {(-1)^n\over n^2(\xi\eta)^{n-1}} [Q_{\mu+1,\nu+1}- &Q_{\mu\nu }-
(Q_{\mu+1,\nu}-Q_{\mu,\nu+1})\cr
&-(Q_{\mu,\nu +1}-Q_{\mu+1,\nu}) + Q_{\mu,\nu}-Q_{\mu+1,\nu+1}] = 0 \cr}
$$
and $\{q_\nu(\xi),q_\mu(\eta)\} = 0$.
 
Turning now to the computation of $\{p_\nu(\xi),q_\mu(\eta)\}$, first note
that the local part vanishes when $\mu\ne \nu$.  As before, this follows from
the block diagonal structure of $a$ or from the structure of the factor
$[sgn(j-l)+sgn(m-k)]$.  When $\mu = \nu$ and $\xi$ and $\eta$ lie on the same
ray,
$$
\eqalign{[\log a_{\mu\mu}(\xi),\log a_{\mu+1,\mu}(\eta)] &= {1\over
a_{\mu\mu}a_{\mu+1,\mu}} [a_{\mu\mu},a_{\mu+1,\mu}]\cr 
&= {-\pi i\over n(-\xi)^{n-1}} \delta (\xi-\eta )sgn(-1).\cr}
$$
We leave it to the reader to verify that the other three terms give the same
result, and that $[\Delta_\mu (\xi),q_\mu(\eta)] = 0$; hence
$$
[p_\nu(\xi),q_\mu(\eta)] = \delta_{\mu \nu}\  \xi \  \delta (\xi - \eta)
= \delta_{\mu \nu} \tilde \delta (\xi/ \eta).
$$ 
Turning now to the nonlocal component, (and neglecting the phase factor
$\alpha^\nu$ in $q_\nu$), we have
$$
<\log {a_{\mu\mu}a_{\mu+1,\mu+1}\over\Delta_\mu} ,\log {a_{\nu+1,\nu}\over
a_{\nu,\nu+1}}>
= {a_{\mu,\mu+1}a_{\mu+1,\mu}\over\Delta_\mu}\times
$$
$$
<\log a_{\mu,\mu+1}a_{\mu+1,\mu} -\log a_{\mu\mu}a_{\mu+1,\mu+1},\log
a_{\nu+1,\nu} -\log a_{\nu,\nu+1}>.
$$
From (2.14)
$$
<\log a_{\mu,\mu+1}a_{\mu+1,\mu} -\log a_{\mu\mu}a_{\mu+1,\mu+1},\log
a_{\nu+1,\nu}>
\eqno(3.2)
$$
$$
= {(-1)^n\over n^2(\xi\eta)^{n-1}}\Bigl[ Q_{\mu,\nu+1} -Q_{\mu+1,\nu+1}-
Q_{\mu,\nu}+Q_{\mu,\nu}- Q_{\mu+1,\nu+1}+Q_{\mu+1,\nu}\Bigr] = 0
$$
and the same holds for the quantity $<\cdot \ ,\log a_{\nu,\nu+1}>$.  Hence 
$<p_\mu(\xi),q_\nu (\eta )> = 0$ and the proof of Theorem 3.1 is
complete.
\vskip .3cm
  
The last computation shows that for all $\mu,~\nu$,
$$
<p_\mu,a_{\nu ,\nu +1}>=<p_\mu, a_{\nu +1,\nu}> = 0.
\eqno(3.3)
$$
This fact will be used in \S5. 
\vskip .3cm 

\noindent{\bf 4. Action-angle variables; self-adjoint case.}  

\vskip .3cm
The canonical coordinates constructed in \S3 are complex; but when $L$ is
Hermitian symmetric, the scattering data possesses an additional symmetry
which allows us to construct real canonical coordinates.  In the following,
$R$ denotes the order reversing matrix
$$
R = \sum^n_{k=1} e_{k,n-k+1}
$$
corresponding to the permutation $(1,n)(2,n-1)\dots$; as before $e_{jk}$ is the
matrix with a 1 in the $jk$ place and zeroes elsewhere.
\vskip .3cm
\noindent{\bf Lemma 4.1}\quad {\it There is a strictly lower triangular
matrix $\Theta (x),\Theta (x)\to 0$ as $\vert x\vert\to\infty$, such that
$$
m_{L^*}(x,z) = (I+\Theta(x))nz^{n-1}R[m_L(x,\bar z)^{-1}]^*RJ(z)^*
\eqno(4.1)
$$
for any $z\in C\backslash\Sigma$ such that $m(x,z)$ exists.}
\vskip .3cm 

In this lemma the two wave functions are given in their local
representations.  When $n$ is even, $\bar z=\alpha^kz$ for some $k$; while if
$n$ is odd, $\bar z = -\alpha^kz$ for some $z$.  Therefore, by the invariance
of the wave functions under the action $z\to\alpha z$ of Lemma 1.2 we may
write (4.1) as
$$
m_{L*}(x,z) = (I+\Theta(x))n\ z^{n-1}R[m_L(x,\sigma^n(z))]^*RJ(z)^* 
\eqno(4.2)
$$ 
where
$$
\sigma (z) = -\bar z.
$$
This expresses the relationship between the wave functions for $L$ and $L^*$
at points in the fundamental domain.  In [BDT] this identity is used to
establish symmetries of the jump data when $L$ is self-adjoint.  We modify
this result here to make it directly applicable to the form of the scattering
data (1.12).  We prove the following:
\vskip .3cm
\noindent{\bf Theorem 4.2.}\quad {\it Let $\xi$ belong to the ray $\Sigma_j$
and $\bar\xi$ to $\Sigma_k$.  Then}
$$
Ra_{L^*}(\bar\xi)R = \pi_jJ^{-1}_j[a_L(\xi)^{-1}]^*J_j\pi_j .
\eqno(4.3)
$$
 
The matrices $a_L$ and $a_{L^*}$ in (4.3) are in their local representations;
$\pi_j$ is the permutation matrix associated with the ray $\Sigma_j$; and
$J_j = \text{diag}(\alpha_j,\dots\alpha_n)$ in the $j^{th}$
ordering.  We first need to prove the analogue of Lemma 4.1 for the wave
functions $m(x,\xi)$ defined for $\xi$ in $\Sigma$.  We indicate the necessary
modification of the proof of Theorem 9.5 in [BDT, pp. 41, 42].  When
$\xi$ belongs to $\Sigma_j$, $m_L(x,\xi)$ is the unique solution of
$$
Dm = J_\xi m - \xi mJ_j + qm,
$$
$$
\lim_{x\to -\infty} m = d(\xi)\Lambda_j , \quad\quad \sup_x\vert m(x,\xi)\vert
<+\infty.
$$
 
\noindent{\bf Lemma 4.3.}\quad {\it Let $\xi\in\ \sigma_j$, and
$\xi\in\Sigma_k$.  Then}
$$
\eqalign{\text{(i)}\ \pi_kJ_k\pi_k &= J_{k+1} ;\ \Lambda_{k+1}\pi_k =
\Lambda_k\cr
\text{(ii)}\ RJ^*_jR &= J_{k+1} ;\ R\Lambda_jR = \Lambda_{k+1}J_{k+1} =
\Lambda_kJ_k\pi_k\cr
\text{(iii)}\ R\pi_jR &= \pi_k ;\ R\Pi_\xi R = \Pi_{\bar\xi} .\cr}
$$
 
\noindent Proof:  The first statement in (i) is a consequence of the defining
relation of the permutations $\pi_j$ and of the ordering of the roots of
unity; the second statement is a result of the fact that only the columns of
$\Lambda$ are permuted in changing representations.  From (1.6),
$$
Re\ i\bar\alpha_n\bar z > Re\ i\bar\alpha_{n-1}\bar z >\dots > Re\
i\bar\alpha_1\bar z \quad\quad \text{for}\ \bar z\in\Omega
$$
and it follows that complex conjugation reverses the ordering of
the roots.  Hence
$$
RJ(z)R = J(\bar z)^*,\quad R\Pi_{\bar\xi} R = \Pi_\xi,\quad \text{and}\ \ R\pi_j R
= \pi_k .
$$
This establishes (iii).  A direct calculation shows that $R\Lambda (z)R =
\Lambda (z)J(z)$.  As $z\to\xi\in\Sigma_j$ from $\Omega_j,\bar
z\to\bar\xi\in\Sigma_k$ from $\Omega_{k+1}$; so the identities in (ii)
follow.  This completes the proof of Lemma 4.3.
 
Repeating the derivation of (9.14, [BDT]) with the wave function \hfill\break
$m(x,\xi)$ we find that
$$
g(x) = n\bar\xi^{n-1} R[m_L(x,\xi)^{-1}]^*RJ_{k+1}^*
$$
satisfies
$$
Dg = J_{\bar\xi} g-\bar\xi gJ_{k+1} +Rq^*Rg,\quad\quad 
g(x)\to\Lambda_{k+1}\quad \text{as}\ x\to -\infty.
$$
Hence $w = g\pi_k$ satisfies
$$
Dw = J_{\bar\xi}w - \bar\xi wJ_k + Rq^*Rw,\quad\quad w\to\Lambda_k\
\text{as}\ x\to -\infty.
$$
As in Theorem 9.5 of [BDT], there exists a matrix function $G(x) = (I+\Theta (x))$
such that
$$
m_{L^*}(x,\bar\xi) = G(x)n\bar\xi^{n-1} R[m_L(x,\xi)^{-1}]^*RJ^*_{k+1}\pi_k .
\eqno(4.5)
$$
\vskip .3cm
\noindent Proof of Theorem 4.2: From (1.12) we have
$$
\eqalign{a_{L^*}(\xi) &= \lim_{x\to\infty} \Pi_{\bar\xi}e^{-ix\bar\xi
J_k}\phi_{L^*}(x,\bar\xi)\cr 
&= \lim_{x\to\infty} \Pi_{\bar\xi} e^{-ix\bar\xi J_k}\Lambda^{-
1}_{\bar\xi}\Phi_{L^*}(x,\bar\xi)\cr}
$$
where
$$
\Phi_{L^*}(x,\bar\xi) = m_{L^*}(x,\bar\xi)e^{ix\bar\xi J_k}\quad\quad
\text{and}\quad\quad \Lambda_{\bar\xi} = d(\bar\xi)\Lambda_k .
$$
From (4.5),
$$
a_{L^*}(\bar\xi) = \lim_{x\to\infty} \Pi_{\bar\xi}e^{ix\bar\xi J_k}
\Lambda^{-1}_{\bar\xi} G(x)n\bar\xi^{n-1} R[m_L(x,\xi)^{-
1}]^*RJ_{k+1}^*\pi_ke^{ix\bar\xi J_k} .
$$
$$
= R \lim_{x\to\infty} \Pi_\xi \pi_jJ^*_j (e^{ix\xi J_j})^*\Lambda_{\xi^*}
[m(x,\xi)^{-1}]^*J_j (e^{-ix\xi J_j})^*\pi_jR
$$
by repeated use of Lemma 4.3.  Since $\pi_j$ and $\Pi_\xi$ have the same
block structure, they commute, so
$$
Ra_{L^*}(\bar\xi)R = \pi_jJ^{-1}_j[\lim_{x\to\infty} \Pi_\xi e^{-ix\xi
J_j}m_L(x,\xi)^{-1}\Lambda_\xi e^{ix\xi J_j}]J_j\pi_j .
$$
Now $\Psi\Phi = a$ and $\Psi = \psi\Lambda^{-1}_\xi$ so $a^{-1}\psi = e^{-
ix\xi J_j}m^{-1}_L\Lambda_\xi$; therefore
$$
\lim_{x\to\infty} \Pi_\xi e^{-ix\xi J_j}\psi = 1
$$
and Theorem 4.2 follows. 
\vskip .3cm

Complex conjugation does not preserve $\Sigma_0 \cup \Sigma_1$, where the 
canonical coordinates were defined in $\S3$, but we may combine it with the
invariance under rotations of Lemma 2.1 to obtain an involution. For
$\xi$ in  $\Sigma_0 \cup \Sigma_1$, let $\xi^*$ be the unique point in
$\Sigma_0 \cup \Sigma_1$ which is on the same orbit as $\xi$ under rotations
by $\alpha$. When $n$ is even, $\xi^*=\xi$; while when n is odd, the involution
interchanges $\Sigma_0$ and $\Sigma_1$. In both cases $\vert \xi^* \vert = 
\vert \xi \vert$.
\vskip .3cm
\noindent{\bf Theorem 4.4.} {\it For $L = L^*$:}
$$
\overline{a_{11}(\xi)} = {1\over a_{nn}(\xi^*)},\quad\quad
\overline{\Delta_\nu(\xi)} = {1\over\Delta_{n-\nu}(\xi^*)},
$$
$$
\eqno(4.6)
$$
$$
\overline{p_\nu(\xi)} = p_{n-\nu}(\xi^*),\quad\quad \overline{q_\nu(\xi)} =
q_{n-\nu}(\xi^*).
$$
\vskip .3cm 
\noindent Proof:  The first identity follows directly from (4.3).  To prove
the remaining identities we use the relation
$$
B_{n-\nu}(\xi^*) = {1\over\overline{\Delta_\nu(\xi)}} \left[ \matrix \bar
a_{\nu+1,\nu+1} & -\alpha^{-\nu}\bar a_{\nu+1,\nu}\\ 
-\alpha^\nu \bar a_{\nu,\nu+1} & \bar a_{\nu\nu} \endmatrix \right] (\xi),
\eqno(4.7)
$$
which is a block by block restatement of (4.3).  In deriving (4.7) recall
that $J_0 =\text{diag}(1,\alpha,\alpha^{-1},\dots )$ and $J_1 =
\text{diag}(\alpha, 1,\alpha^2,\alpha^{-1},\dots )$; we work in the $0^{th}$
representation.
 
The second identity in (4.6) is immediate from (4.7).  We leave it to the
reader to verify the third; while for that involving the $q_\nu$'s, we have
$$
q_{n-\nu}(\xi^*) = {1\over 4\pi i} \log {a_{n-\nu+1,n-\nu}\over a_{n-\nu,n-
\nu+1}} \alpha^{n-\nu}
= {1\over 4\pi i}\log \alpha^\nu{\bar a_{\nu,\nu+1}\over\alpha^{-\nu}\bar
a_{\nu+1,\nu}}\alpha^{n-\nu} 
$$
$$
= {1\over 4\pi i}\log \alpha^\nu {\bar
a_{\nu,\nu+1}\over\bar a_{\nu+1,\nu}}
= -{1\over 4\pi i}\log\alpha^{-\nu}{\bar a_{\nu+1,\nu}\over\bar
a_{\nu,\nu+1}} = \bar q_\nu (\xi).
$$
\vskip .3cm
 
These symmetries allow us to construct real canonical variables in the self-
adjoint case.  For odd n they are
$$
\eqalign{\sqrt 2 Re\ p_\nu(\xi),\quad &\sqrt 2 Re\ q_\nu(\xi),\cr
\sqrt 2 Im\ p_\nu(\xi),\quad &\sqrt 2Im\ q_\nu(\xi),\quad\quad 1\leq \nu <\dots
<n/2.\cr}
$$
When $n=2m$ is even, $p_m$ and $q_m$ are real and 
these variables are to be added to the above list. 
These real variables also satisfy the canonical commutation relations when
$\xi$ and $\eta$ lie on the rays $\Sigma_0$ and $\Sigma_1$.  For example
when $\mu$, $\nu<n/2$,
$$
\{\sqrt 2\ Re\ p_\nu(\xi),\sqrt 2\ Re\ q_\mu(\eta)\} = \{{p_\nu(\xi)+p_{n-
\nu}(\xi^*)\over\sqrt 2},\ {q_\mu(\eta)+q_{n-\mu}(\eta^*)\over\sqrt 2}\}
$$
$$
={1\over 2}[\delta_{\mu\nu}\tilde \delta (\xi/\eta) + \delta_{\nu,n-\mu}
\tilde \delta(\xi/
\eta^*)+\delta_{n-\nu,\mu}\tilde \delta (\xi^*/\eta)+\delta_{n-\nu,n-\mu}
\tilde \delta (\xi^*/ \eta^*)]
$$
$$
=\delta_{\nu\mu}\tilde \delta (\xi/ \eta)+\delta_{\nu,n-\mu}
\tilde \delta(\xi/ \eta^*) =
\delta_{\nu\mu}\tilde \delta (\xi/ \eta)
$$
since $\delta_{\nu,n-\mu}=0$ when $\mu,\nu<n/2$.  

\vskip .3cm
\noindent{\bf 5. Complete Integrability of the Gel'fand-Dikii Flows.}
\vskip .3cm
We may now demonstrate the complete integrability of the Gel'fand-Dikii flows. 
Lax pairs for the flows (1.1) can be obtained [BDT] as a pair of first order
operators $D_x = D-zJ-q_z$ and $D_{k,t} = D_t-(zJ)^k-B_{n,k}$, where $D_t = -
id/dt$, and the $B_{n,k}$ are $n\times n$ matrices whose entries are
polynomials in $z,u_0,\dots u_{n-1}$ and their derivatives up to order $k-
1$.  Under the flows given by the zero curvature condition $[D_x,D_{k,t}] =
0$, the scattering matrix $a$ evolves according to
$$
a(\xi,t) = e^{i(\xi J)^kt}a(\xi,0) e^{-i(\xi J)^kt}.
\eqno(5.1)
$$
As in the case of integrable systems built on first order $n\times n$
isospectral operators, the scattering transform {\it decouples} the modes, in
the sense that the evolution of the scattering matrix at $\xi\in\Sigma$ is
independent of $\xi'\ne\xi$. It is easily seen that the pointwise flows (5.1)
 are Hamiltonian with
respect to the pointwise Poisson bracket
$$
(a_{jk},a_{lm}) = { \pi i\over n(-\xi)^n } a_{lk}a_{jm}[sgn(l-j)+sgn(k-m)]. 
\eqno(5.2)
$$
Relative to the bracket (5.2), the local component of (2.14) is
$(a_{jk}(\xi),\ a_{lm}(\eta))\tilde \delta (\xi/\eta)$, since 
$\xi~ \delta(\xi - \eta)=\tilde \delta(\xi/\eta).$
 
In what sense do the flows (5.1) comprise a ``complete set of commuting
Hamiltonian flows?''  For fixed $\xi\in\Sigma$ the minimal scattering data at
$\xi$ is given by the upper and lower triangular factors $v^\nu_\pm$ of
$B_\nu$ with diag $v^\nu_\pm = 1$ (cf. [BC] or [BS1]).  Each $v^\nu_\pm$ has
one nontrivial entry, so the scattering data on $\Sigma$ consists of a total
of $2(n-1)$ functions defined on the two half lines $\Sigma_0$ and
$\Sigma_1$.  (Compare with the $n-1$ functions $u_0,\dots u_{n-2}$ defined on
the real line).  There are $[n/2]$ degrees of freedom for each
$\xi\in\Sigma_1$ and $n-[n/2]-1$ degrees of freedom at each point
$\xi\in\Sigma_0$.  The question of complete integrability of the 
Gel'fand-Dikii flows is now reduced to a finite dimensional problem, namely that of
finding $n-1$ pointwise independent commuting flows on the scattering data
for each fixed $\xi\in\Sigma$.  It is clear that, for fixed $\xi\in\Sigma$,
the flows (5.1) comprise just such a maximal set of commuting flows.
 
The Hamiltonians for the flows (5.1) are generated by a linear combination of
the $p_\nu$'s.  Straightforward computations using the properties of the
brackets (e.g. the derivation property) show
$$
\eqalign{(p_\nu,a_{\mu,\mu+1}) &= -2\pi i\delta_{\nu\mu}a_{\mu,\mu+1},\cr
(p_\nu,a_{\mu+1,\mu}) &= 2\pi i\delta_{\nu\mu}a_{\mu+1,\mu},\quad\quad
(p_\nu,a_{\mu\mu}) = 0.\cr}
$$
Therefore the flow defined by

$$
D_ta_{jk} = (a_{jk},{1\over 2\pi i}p_\nu),\quad D_t = {1\over i} {d\over dt}
$$
is given by $B_\nu(t) = \exp\{i\sigma_\nu t\}B_\nu(0)\exp\{-i\sigma_\nu t\}$
and $B_\mu(t) = B_\mu(0)$ if $\mu\ne\nu$, where $B_\nu$ is the
$2\times 2$ diagonal block and $\sigma_\nu$ is the diagonal matrix with a 1/2
in the $\nu^{th}$ entry and -1/2 in the $(\nu+1)^{st}$ entry.
 
For $\xi$ in $\Sigma_j$, let $av_{\pm}=a_{\pm}$ (assuming this 
factorization exists), where $v_+,~a_-$ $(v_-,~a_+)$ are lower (upper) triangular
in the $j^{th}$ representation; and diag $v_{\pm}=1$. 
The matrix $v(\xi)$ introduced in \S1 as the scattering data equals 
$v_-^{-1}v_+$ ([BC]). Define $\delta_{\pm}=$diag $a_{\pm}.$
By direct computation (cf. Lemma 4.2 of [BS1]) we obtain
$$
\delta^{-1}_-\delta_+ = \text{diag}\Bigl( {\Delta_1\over
a_{11}a_{22}},{a_{11}a_{22}\over\Delta_1},{\Delta_2\over a_{33}a_{44}},\dots
\Bigr)
$$
on the odd numbered rays, with a corresponding result holding on
the even numbered rays. By (3.1)
$$
-{n(-\xi)^n \over 2\pi i} \log\delta^{-1}_-\delta_+ = {1\over 2\pi i}
\text{diag}(p_1,-p_1,p_3,-p_3,\dots ) 
\eqno(5.3)
$$
on the odd rays, again with a corresponding result on the even rays.  
Therefore the pointwise flows (5.1) are generated by
$$
-{n(-\xi)^n \over 2\pi i} \text{tr}[(\xi J)^k\log\delta^{-1}_-\delta_+] .
$$
The flows (5.1) are generated by the direct integral of these pointwise
Hamiltonians
$$ 
\eqalign{ {\Cal H}_k &= {n \over 2\pi i} \int_{\Sigma_0\cup \Sigma_1}
(-\xi)^n \text{tr}[(\xi J)^k\log\delta^{-1}_-
\delta_+]{d\xi \over \xi}\cr
&= {1\over 2\pi i}\int_{\Sigma}(-\xi)^n \text{tr}[\xi J)^k\log\delta^{-1}_-
\delta_+]{d\xi \over \xi}.
\cr} 
\eqno(5.4)
$$  
Indeed, consider the functionals
$$
F_{lm}=\int_\Sigma a_{lm}(\xi,t)f(\xi){d\xi \over \xi}
$$
and the associated flows defined by 
$D_t F_{lm}=\{F_{lm},{\Cal H}_k\}.$ By (3.3) the nonlocal part of these brackets vanishes,
and by (5.2) and (5.3),
$$
\{F_{lm},{\Cal H}_k\} = \int_\Sigma f(\xi)a_{lm}(\xi,t)(\alpha_l^k-\alpha_m^k)\xi^k
{d\xi \over \xi}.
$$    
Since f is an arbitrary element of $\Cal D (\Sigma)$,
$$
D_ta_{lm}=(\alpha_l^k-\alpha_m^k)\xi^ka_{lm},
$$
and the corresponding flows are precisely the Gel'fand-Dikii flows (5.1).

These integrals
are well defined if the potential $q$ is in the Schwartz class; for then ([BDT], 
Corollary 5.18) $\log\delta^{-1}_-\delta_+$ tends to zero rapidly as $\xi$ tends to infinity
along the rays of $\Sigma$. Moreover ([BC]) the diagonal matrices $\delta_{\pm}$
are the boundary values on $\Sigma$ of the
sectionally meromorphic function $\delta (z)$ obtained as
$$
\delta (z) = \lim_{x\to \infty} \phi (x,z)e^{-ixzJ(z)}.
$$
In the absence of discrete components of the scattering data, $\delta (z)$
is holomorphic in $\Omega_{\pm}$ and its entries have no zeroes.
 
For integrable systems based on $n\times n$ first order isospectral operators
an $(n-1)$-parameter family of generating functions for the Hamiltonians 
is given by tr$[\mu\log\delta(z)]$, $\mu$ a traceless diagonal matrix  
([S],[BC]).  We have a similar result here, with a single generator.
\vskip .3cm
\noindent{\bf Theorem 5.1}\quad {\it The generating function for the
Gel'fand-Dikii flows is}
$$
{\Cal H}(z) = (-1)^{n+1}z^n\text{tr}[(I+J+J^2 +\dots +J^{n-1})\log\delta (z)].
$$
\vskip .2cm 
Note that $n^{-1}(I+J+\dots +J^{n-1})$ projects onto the null space of $J - I$. 
The generating function has branch points at the zeroes
of $\delta_j (z)$, which correspond to the discrete part of the
scattering data; but we again omit this case and prove Theorem 5.1 under the
assumption that the scattering data has no such discrete part.  Under
those assumptions ${\Cal H}(z)$ is single valued and analytic in each
component of $C\backslash\Sigma$.
\vskip .3cm
\noindent Proof:\quad The sectionally analytic function $\Phi_k(z) = \text{tr}
[J^k(z)\log\delta (z)]$ is $O(z^{-1})$ as $z\to\infty$.  Its
jump across $\Sigma$ is $[\Phi_k] = \text{tr}[ J(\xi)\log\delta^{-1}_-\delta_+]$
(using the same representation for $\delta_-$ and $\delta_+$ on $\Sigma$) so
by the Plemelj formulas
$$
\Phi_k(z) = {1\over 2\pi i} \int_\Sigma {\text{tr}[J^k(\xi)\log\delta^{-1}_-
\delta_+]\over\xi -z} d\xi  .
$$
(There is a weak singularity at $\xi = 0$: $\delta (z)$ behaves
like a power of $z$ as $z$ tends to zero in the closure of any of the sectors
$\Omega_j$ [BDT]; so $\log\delta^{-1}_-\delta_+$ has a logarithmic
singularity at the origin.)  Since $zJ(z)$ is invariant under $z\to\alpha z$,
it follows that $\Phi_k(\alpha z) = \alpha^{-k}\Phi_k(z)$, and
$$
n\Phi_k(z) = {1\over 2\pi i} \int_\Sigma \sum^{n-1}_{j=0}
{\alpha^{jk}\over\xi -\alpha^jz} \text{tr}[ J^k(\xi)\log\delta^{-1}_-\delta_+] d\xi.
$$
From the identity
$$
\sum^{n-1}_{j=0} {\alpha^{jk}\over\xi-\alpha^jz} = n{\xi^{k-1}z^{n-
k}\over\xi^n-z^n}
$$
which may be verified by comparing the residues of both sides with respect to
the variable $\xi$, it follows that
$$
\text{tr}[J^k\log\delta(z)] = -{1\over z^k}{1\over 2\pi i}\int_\Sigma {\xi^k \over 1-(\xi
/z)^n} \text{tr}[ J^k\log\delta^{-1}_-\delta_+ ]{d\xi\over \xi}
$$
$$
\eqalign{& \sim -\sum^\infty_{j=0} {1\over z^{nj+k}} {1\over 2\pi i} \int_\Sigma
\xi^{nj+k}\text{tr}[ J^k\log\delta^{-1}_-\delta_+] {d\xi\over \xi} \cr
&= (-1)^{n+1}\sum^\infty_{j=0} {1\over z^{nj+k}} {1\over 2\pi i}\int_\Sigma
(-\xi )^n \text{tr}[(\xi J)^{n(j-1)+k}\log\delta^{-1}_-\delta_+]{ d\xi \over \xi}\cr
&= (-1)^{n+1}z^{-n} \sum^\infty_{j=-1} {1 \over z^{nj+k}} {\Cal H}_{nj+k} .\cr}
$$
The result is obtained by summing from $k=0$ to $n-1$.
\vskip .3cm
We would like to thank Boris Konopelchenko for helpful discussions during the preparation
of this paper.
 
\enddocument